\documentclass[12pt]{article}%
\usepackage{graphicx}
\usepackage[intlimits]{amsmath}
\usepackage{latexsym}
\usepackage{amsfonts}
\usepackage{amssymb}%
\makeatletter
\newfont{\footsc}{cmcsc10 at 8truept}
\newfont{\footbf}{cmbx10 at 8truept}
\newfont{\footrm}{cmr10 at 10truept}
\newtheorem{theorem}{Theorem}

\newtheorem{lemma}[theorem]{Lemma}

\newtheorem{proposition}[theorem]{Proposition}

\newenvironment{proof}[1][Proof]{\noindent{\textbf {#1}  }}  {\hfill$\Box$\bigskip}

\begin{document}

\title{Numerical radius and zero pattern of matrices}
\author{Vladimir Nikiforov\\{\small Department of Mathematical Sciences, University of Memphis, }\\{\small Memphis TN 38152, USA, email: }\textit{vnkifrv@memphis.edu}}
\maketitle

\begin{abstract}
Let $A$ be an $n\times n$ complex matrix and $r$ be the maximum size of its
principal submatrices with no off-diagonal zero entries. Suppose $A$ has zero
main diagonal and $\mathbf{x}$ is a unit $n$-vector. Then, letting $\left\Vert
A\right\Vert $ be the Frobenius norm of $A,$ we show that
\[
\left\vert \left\langle A\mathbf{x},\mathbf{x}\right\rangle \right\vert
^{2}\leq\left(  1-1/2r-1/2n\right)  \left\Vert A\right\Vert ^{2}.
\]
This inequality is tight within an additive term $O\left(  n^{-2}\right)  .$

If the matrix $A$ is Hermitian, then%
\[
\left\vert \left\langle A\mathbf{x},\mathbf{x}\right\rangle \right\vert
^{2}\leq\left(  1-1/r\right)  \left\Vert A\right\Vert ^{2}.
\]
This inequality is sharp; moreover, it implies the Tur\'{a}n theorem for
graphs. \bigskip


\bigskip

\end{abstract}

\section{Introduction}

Let $G$ be a simple graph, $\mu\left(  G\right)  $ be the spectral radius of
its adjacency matrix, $\omega\left(  G\right)  $ be the maximum size of its
complete subgraphs, and $e\left(  G\right)  $ be the number of its edges. In
\cite{Nik02} it is shown that
\begin{equation}
\mu^{2}\left(  G\right)  \leq\left(  2-\frac{2}{\omega\left(  G\right)
}\right)  e\left(  G\right)  . \label{mainin}%
\end{equation}
The aim of this note is to extend this result to square matrices with zero
main diagonal.

Let $\eta\left(  A\right)  $ be the numerical radius of a square matrix $A,$
i.e.,
\[
\eta\left(  A\right)  =\max_{\left\Vert \mathbf{x}\right\Vert =1}\left\vert
\left\langle A\mathbf{x},\mathbf{x}\right\rangle \right\vert .
\]
The value $\eta\left(  A\right)  $ has been extensively studied, see, e.g.,
\cite{Fua03}-\cite{MeKu05}, \cite{MPT02} and their references.

Given a complex matrix $A=\left\{  a_{ij}\right\}  ,$ write $\left\Vert
A\right\Vert $ for its Frobenius's norm, i.e., $\left\Vert A\right\Vert
=\sqrt{\sum_{i,j}\left\vert a_{ij}\right\vert ^{2}}.$ We are interested in
upper bounds on $\eta\left(  A\right)  $ in terms of $\left\Vert A\right\Vert
.$ It is easy to see that $\eta\left(  A\right)  \leq\left\Vert A\right\Vert $
with equality holding, e.g., if $A$ is a constant matrix. In this note we give
conditions for the zero pattern of a square matrix $A$ that imply $\eta\left(
A\right)  \leq\left(  1-c\right)  \left\Vert A\right\Vert $ for some
$c\in\left(  0,1\right)  $ independent of the order of $A.$

Given a square matrix $A,$ let $\omega\left(  A\right)  $ be the maximum size
of its principal submatrices with no off-diagonal zero entries.

Note that if $A$ is the adjacency matrix of a graph $G,$ then $\omega\left(
A\right)  =\omega\left(  G\right)  ,$ $\mu\left(  G\right)  =\eta\left(
A\right)  ,$ and $\left\Vert A\right\Vert ^{2}=2e\left(  G\right)  .$ Thus,
the following theorem extends inequality (\ref{mainin}).

\begin{theorem}
\label{th1} For every Hermitian matrix $A$ with zero main diagonal,
\begin{equation}
\eta^{2}\left(  A\right)  \leq\left(  1-\frac{1}{\omega\left(  A\right)
}\right)  \left\Vert A\right\Vert ^{2}. \label{in1}%
\end{equation}
Inequality (\ref{in1}) is sharp: for all $n\geq r\geq2,$ there exists an
$n\times n$ symmetric $\left(  0,1\right)  $-matrix $A$ with zero main
diagonal and $\omega(A)=r$ such that equality holds in (\ref{in1})$.$
\end{theorem}

Note that inequality (\ref{in1}) implies a concise form of the fundamental
theorem of Tur\'{a}n in extremal graph theory (see \cite{Bol98} for details).
Indeed, if $A$ is the adjacency matrix of a graph $G$ with $n$ vertices and
$m$ edges, then
\[
\left(  2m/n\right)  ^{2}\leq\eta^{2}\left(  A\right)  \leq\left(
2-2/\omega\left(  A\right)  \right)  m=2\left(  1-1/\omega\left(  G\right)
\right)  m,
\]
and so,
\begin{equation}
m\leq\left(  1-\frac{1}{\omega\left(  G\right)  }\right)  \frac{n^{2}}{2}.
\label{maxmu1}%
\end{equation}
Moreover, inequality (\ref{in1}) follows from a result of Motzkin and Straus
\cite{MoSt65}, following in turn from (\ref{maxmu1}) (see \cite{Nik06a} for
details). The implications
\[
\left(  \ref{in1}\right)  \Longrightarrow\left(  \ref{maxmu1}\right)
\Longrightarrow\text{MS}\Longrightarrow\left(  \ref{in1}\right)
\]
justify regarding inequality (\ref{in1}) as a matrix form of Tur\'{a}n's theorem.

We state without a proof a characterization of Hermitian matrices for which
equality holds in (\ref{in1}).

\begin{proposition}
Let $A=\left\{  a_{ij}\right\}  $ be an $n\times n$ Hermitian matrix with zero
main diagonal with $\omega\left(  A\right)  =r\geq2.$ Then the equality
$\eta^{2}\left(  A\right)  =\left(  1-1/r\right)  \left\Vert A\right\Vert
^{2}$ holds if and only if there exist a complex number $c\neq0,$ a partition
$\left[  n\right]  =\cup_{i=0}^{r}N_{i},$ and a unit vector $\mathbf{x}%
=\left(  x_{1},\ldots,x_{n}\right)  $ such that:

(i) $x_{i}=0$ for all $i\in N_{0}.$

(ii) $\sum_{i\in N_{i}}\left\vert x_{i}\right\vert ^{2}=1/r$ for all $1\leq
i\leq r.$

(iii) $a_{ij}=cx_{i}\overline{x}_{j}$ for all $1\leq i<j\leq n.$
\end{proposition}

It turns out that Theorem \ref{th1} has analogues for non-Hermitian matrices
as well.

\begin{theorem}
\label{th2} For every complex $n\times n$ matrix $A$ with zero main diagonal,
\begin{equation}
\eta^{2}\left(  A\right)  \leq\left(  1-\frac{1}{2\omega\left(  A\right)
}-\frac{1}{2n}\right)  \left\Vert A\right\Vert ^{2}. \label{in2}%
\end{equation}
Inequality (\ref{in2}) is tight: for all $n\geq r\geq2,$ there exists an
$n\times n$ matrix $A$ with zero main diagonal and $\omega(A)=r$ such that
\[
\eta^{2}\left(  A\right)  \geq\left(  1-\frac{1}{2\omega\left(  A\right)
}-\frac{1}{2n}+O\left(  n^{-2}\right)  \right)  \left\Vert A\right\Vert ^{2}.
\]

\end{theorem}

Let $P_{n}$ be the set of vectors $(x_{1},\ldots,x_{n})$ with $x_{1}%
\geq0,\ldots,x_{n}\geq0,$ and $x_{1}+\cdots+x_{n}=1.$ Recall a result of
Motzkin and Straus \cite{MoSt65}: if $A$ is the adjacency matrix of a graph
$G$ of order $n,$ and $\mathbf{x}\in P_{n},$ then%
\begin{equation}
\left\langle A\mathbf{x},\mathbf{x}\right\rangle \leq1-1/\omega\left(
G\right)  . \label{MS_in}%
\end{equation}

We shall need the following extension of this result.

\begin{lemma}
\label{cl1} For every square $\left(  0,1\right)  $-matrix $A$ of size $n$
with zero main diagonal and every $\mathbf{x}\in P_{n},$%
\begin{equation}
\left\langle A\mathbf{x},\mathbf{x}\right\rangle \leq1-\frac{1}{2\omega\left(
A\right)  }-\frac{1}{2n}. \label{MS_in1}%
\end{equation}
Inequality (\ref{MS_in1}) is tight: for all $n\geq r\geq2,$ there exists a
square $\left(  0,1\right)  $-matrix $A$ of size $n$ with zero main diagonal
and $\omega\left(  A\right)  =r$ such that,
\[
\left\langle A\mathbf{x},\mathbf{x}\right\rangle =1-\frac{1}{2r}-\frac{1}%
{2n}+O\left(  n^{-2}\right)
\]
for some $\mathbf{x}\in P_{n}.$
\end{lemma}

\section{Proofs}

\begin{proof}
[\textbf{Proof of Lemma \ref{cl1}}]Define the $n\times n$ matrix $B=\left\{
b_{ij}\right\}  $ setting $b_{ij}=a_{ij}a_{ji}$ for all $i,j\in\left[
n\right]  ;$ let $C=A-B.$ Note that for every two distinct $i,j\in\left[
n\right]  ,$ we have%
\[
c_{ij}+c_{ji}=a_{ij}+a_{ij}-2a_{ij}a_{ji}\leq1.
\]
We may and shall assume that $c_{ij}+c_{ji}=1$ for all distinct $i,j\in\left[
n\right]  $ with $b_{ij}=0,$ since otherwise some off-diagonal zero entry of
$A$ can be changed to $1$ so that $\omega\left(  A\right)  $ remains the same
and the left-hand side of (\ref{MS_in1}) does not decrease. Hence, for every
$\mathbf{x}=(x_{1},\ldots,x_{n}),$
\[
\left\langle B\mathbf{x},\mathbf{x}\right\rangle +2\left\langle C\mathbf{x}%
,\mathbf{x}\right\rangle =1-\left\Vert \mathbf{x}\right\Vert ^{2}.
\]
Since $B$ is a symmetric $\left(  0,1\right)  $-matrix with zero main
diagonal, the result of Motzkin and Straus implies that
\[
\left\langle B\mathbf{x},\mathbf{x}\right\rangle \leq1-1/\omega\left(
B\right)
\]
for every $\mathbf{x}\in P_{n}.$ Since $\omega\left(  B\right)  =\omega\left(
A\right)  ,$ we find that
\[
\left\langle A\mathbf{x},\mathbf{x}\right\rangle =\left\langle B\mathbf{x}%
,\mathbf{x}\right\rangle +\left\langle C\mathbf{x},\mathbf{x}\right\rangle
=\frac{1}{2}\left(  1-\left\Vert \mathbf{x}\right\Vert ^{2}\right)  +\frac
{1}{2}\left\langle B\mathbf{x},\mathbf{x}\right\rangle \leq1-\frac{1}%
{2\omega\left(  A\right)  }-\frac{1}{2n}.
\]
completing the proof of (\ref{MS_in1}).

Let $G$ be a complete $r$-partite graph whose vertex classes differ in size by
at most $1$. Let $T=\left\{  t_{ij}\right\}  $ be the adjacency matrix of $G;$
set $t_{ij}=1$ for $i<j$ and write $A$ for the resulting matrix. We have
\[
\left\Vert A\right\Vert ^{2}=\binom{n}{2}+\frac{1}{2}\left\Vert T\right\Vert
^{2}=\binom{n}{2}+\binom{r}{2}\frac{n^{2}-\nu^{2}}{r^{2}}+\binom{\nu}{2},
\]
Letting $\mathbf{x}$ to be the $n$-vector $\left(  1/n,\ldots,1/n\right)  \in
P_{n},$ we find that
\begin{align*}
\left\langle A\mathbf{x},\mathbf{x}\right\rangle  &  =\frac{1}{n^{2}%
}\left\Vert A\right\Vert ^{2}=\frac{1}{n^{2}}\left(  \binom{n}{2}+\binom{r}%
{2}\frac{n^{2}-\nu^{2}}{r^{2}}+\binom{\nu}{2}\right) \\
&  =1-\frac{1}{2r}-\frac{1}{2n}+\left(  \frac{\nu^{2}}{2r}-\frac{\nu}%
{2}\right)  \frac{1}{n^{2}}\geq1-\frac{1}{2r}-\frac{1}{2n}-\frac{r}{8n^{2}},
\end{align*}
completing the proof of the lemma.
\end{proof}

\bigskip

\begin{proof}
[\textbf{Proof of Theorem \ref{th1}}]Select $\mathbf{y}=(y_{1},\ldots,y_{n})$
with $\left\Vert \mathbf{y}\right\Vert =1$ and $\eta\left(  A\right)
=\left\vert \left\langle A\mathbf{y},\mathbf{y}\right\rangle \right\vert .$ We
have, by the Cauchy-Schwarz inequality,
\[
\eta^{2}\left(  A\right)  =\left\vert
{\textstyle\sum\limits_{i,j}}
a_{ij}y_{i}\overline{y}_{j}\right\vert ^{2}\leq%
{\textstyle\sum\limits_{i,j}}
\left\vert a_{ij}\right\vert ^{2}%
{\textstyle\sum\limits_{a_{ij}\neq0}}
\left\vert y_{i}\right\vert ^{2}\left\vert y_{j}\right\vert ^{2}=\left\Vert
A\right\Vert ^{2}%
{\textstyle\sum\limits_{a_{ij}\neq0}}
\left\vert y_{i}\right\vert ^{2}\left\vert y_{j}\right\vert ^{2}.
\]
Define a graph $G$ with $V\left(  G\right)  =\left[  n\right]  ,$ joining $i$
and $j$ if $a_{ij}\neq0.$ Obviously, $\omega\left(  G\right)  =\omega\left(
A\right)  $. Since $\left\Vert \mathbf{y}\right\Vert =1,$ the result of
Motzkin and Straus implies that
\[%
{\textstyle\sum\limits_{a_{ij}\neq0}}
\left\vert y_{i}\right\vert ^{2}\left\vert y_{j}\right\vert ^{2}=%
{\textstyle\sum\limits_{ij\in E\left(  G\right)  }}
\left\vert y_{i}\right\vert ^{2}\left\vert y_{j}\right\vert ^{2}\leq1-\frac
{1}{\omega\left(  A\right)  },
\]
completing the proof of (\ref{in1}).

Let $A$ be the adjacency matrix of the union of a complete graph on $r$
vertices and $n-r$ isolated vertices. Since $\omega\left(  A\right)  =r,$
$\eta\left(  A\right)  =r-1,$ and $\left\Vert A\right\Vert ^{2}=r\left(
r-1\right)  ,$ we see that
\[
\eta^{2}\left(  A\right)  =\mu^{2}\left(  A\right)  =\left(  1-1/\omega\left(
A\right)  \right)  \left\Vert A\right\Vert ^{2},
\]
completing the proof of the theorem.
\end{proof}

\bigskip

\begin{proof}
[\textbf{Proof of Theorem \ref{th2}}]Select $\mathbf{y}=(y_{1},\ldots,y_{n})$
with $\left\Vert \mathbf{y}\right\Vert =1$ and $\eta\left(  A\right)
=\left\vert \left\langle A\mathbf{y},\mathbf{y}\right\rangle \right\vert .$
Lemma \ref{cl1} implies that
\begin{align*}
\eta^{2}\left(  A\right)   &  =\left\vert
{\textstyle\sum\limits_{i,j}}
a_{ij}y_{i}\overline{y}_{j}\right\vert ^{2}\leq%
{\textstyle\sum\limits_{i,j}}
\left\vert a_{ij}\right\vert ^{2}%
{\textstyle\sum\limits_{a_{ij}\neq0}}
\left\vert y_{i}\right\vert ^{2}\left\vert y_{j}\right\vert ^{2}=\left\Vert
A\right\Vert ^{2}%
{\textstyle\sum\limits_{a_{ij}\neq0}}
\left\vert y_{i}\right\vert ^{2}\left\vert y_{j}\right\vert ^{2}\\
&  \leq\left(  1-\frac{1}{2\omega(A)}-\frac{1}{2n}\right)  \left\Vert
A\right\Vert ^{2},
\end{align*}
proving (\ref{in2}).

To complete the proof, select $A$ as in the proof of Lemma \ref{cl1}. Hence,
letting $\nu$ be the remainder of $n$ modulo $r,$ we have
\[
\left\Vert A\right\Vert ^{2}=%
{\textstyle\sum\limits_{i,j}}
a_{ij}=\left(  \binom{n}{2}+\binom{r}{2}\frac{n^{2}-\nu^{2}}{r^{2}}+\binom
{\nu}{2}\right)  .
\]
Selecting $\mathbf{x}$ to be the $n$-vector $\left(  n^{-1/2},\ldots
,n^{-1/2}\right)  ,$ as in the proof of Lemma \ref{cl1}, we find that
\[
\eta^{2}\left(  A\right)  \geq\frac{1}{n^{2}}\left\Vert A\right\Vert
^{2}=1-\frac{1}{2r}-\frac{1}{2n}+\left(  \frac{\nu^{2}}{2r}-\frac{\nu}%
{2}\right)  \frac{1}{n^{2}}\geq1-\frac{1}{2r}-\frac{1}{2n}-\frac{r}{8n^{2}},
\]
completing the proof of the theorem.
\end{proof}

\bigskip

\textbf{Concluding remarks}

- The example constructed in the proof of Lemma \ref{cl1} shows that equality
may hold in (\ref{in2}) and (\ref{MS_in1}) whenever $n$ is a multiple of $r$.

- It would be interesting to drop the requirement for zero main diagonal in
Theorem \ref{th1} and \ref{th2}. Note that inequalities (\ref{MS_in}) and
(\ref{MS_in1}) are no longer valid if ones are present on the main diagonal of
$A.$

- Since the spectral radius of a square matrix does not exceed its numerical
radius, Theorem \ref{th1} and \ref{th2} provide upper bounds on the spectral
radius as well.

\end{document}